\documentclass{amsart}
\usepackage{amsmath,amssymb,amscd,latexsym}

\newcommand{\C}{{\mathbb{C}}}

\newcommand{\N}{{\mathbb{N}}}

\newcommand{\Bh}{{\mathcal B}}
\newcommand{\Ch}{{\mathcal C}}

\newcommand{\Fh}{{\mathcal F}}

\newcommand{\Oh}{{\mathcal O}}

\newcommand{\Zh}{{\mathcal Z}}

\newcommand{\Aff}{\mathrm{Aff}}

\newcommand{\be}{\mathbf{1}}

\newcommand{\dist}{\mathrm{dist}}
\newcommand{\dr}{\mathrm{dr}\,}

\newcommand{\halb}{\frac{1}{2}}

\newcommand{\id}{\mathrm{id}}
\newcommand{\ord}{\mathrm{ord}\,}

\newcounter{number}[section]

\newenvironment{nummer}{\refstepcounter{number}{\noindent\arabic{section}.\arabic{number}}}{}

\newcommand{\bn}{\noindent \begin{nummer} \rm}
\newcommand{\en}{\end{nummer}}

\newenvironment{ntheorem}{\noindent {\sc Theorem:} \it}{}
\newenvironment{nlemma}{\noindent {\sc Lemma:} \it}{}
\newenvironment{nprop}{\noindent {\sc Proposition:} \it}{}
\newenvironment{ndefn}{\noindent {\sc Definition:} \it}{}
\newenvironment{ncor}{\noindent {\sc Corollary:} \it}{}

\newenvironment{nproof}{\noindent {\sc Proof:}}{\mbox{}\hfill 
\rule[-.2ex]{.25em}{1.8ex}}

\begin{document}

\title[{Classification, $\Zh$-stability and decomposition rank}]{{On the classification of simple $\Zh$-stable  $C^{*}$-Algebras with real rank zero and finite decomposition rank}}

\author{Wilhelm Winter}
\address{Mathematisches Institut der Universit\"at M\"unster\\
Einsteinstr.\ 62\\ D-48149 M\"unster}

\email{wwinter@math.uni-muenster.de}

\date{January 2006}

\thanks{{\it Supported by:} EU-Network Quantum Spaces - Noncommutative 
Geometry (Contract No. \\
\indent HPRN-CT-2002-00280) and Deutsche 
Forschungsgemeinschaft (SFB 478)}

\setcounter{section}{-1}

\maketitle

\begin{abstract}
We show that, if $A$ is a separable simple unital $C^{*}$-algebra which absorbs the Jiang--Su algebra $\Zh$ tensorially and which has real rank zero and finite decomposition rank, then $A$ is tracially AF in the sense of Lin, without any restriction on the tracial state space. As a consequence, the Elliott conjecture is true for the class of $C^{*}$-algebras as above which, additionally, satisfy the Universal Coefficients Theorem. In particular, such algebras are completely determined by their ordered $K$-theory. They are approximately homogeneous of topological dimension less than or equal to 3, approximately subhomogeneous of topological dimension at most 2 and their decomposition rank also is no greater than 2. 
\end{abstract}

\section{Introduction}

It is the aim of the Elliott classification program to find complete $K$-theoretic invariants for separable simple nuclear $C^{*}$-algebras. For purely infinite $C^{*}$-algebras, this task was accomplished by Kirchberg and Phillips, and there are numerous classification results for special inductive limits of (sub)homogeneous $C^{*}$-algebras, cf.\ \cite{R1} for an overview. However, Toms (in the stably finite case) and R{\o}rdam (in the infinite case) have given examples which show that the Elliott invariant, at least in its current form (cf.\ \cite{R1}, Section 2.2), is not complete for the class of all simple nuclear $C^{*}$-algebras. The construction of these counterexamples is based on techniques introduced by Villadsen; they involve limits of homogeneous $C^{*}$-algebras with `fast dimension growth', i.e., the covering dimension of the underlying base spaces grows faster than the vector space dimension of the matrix algebras which occur as fibers.    

It is interesting to ask for a subclass of separable simple nuclear $C^{*}$-algebras which can be described in a natural way and, at the same time, contains the classes which have been classified so far, but excludes the known counterexamples to the Elliott conjecture. There is growing body of evidence that a very useful notion of `good' behavior for simple $C^{*}$-algebras is the concept of $\Zh$-stability, i.e., the property of absorbing the Jiang--Su algebra $\Zh$ tensorially. This algebra was constructed by Jiang and Su in \cite{JS1}; it is a simple, stably finite and infinite dimensional $C^{*}$-algebra which is $KK$-equivalent to the complex numbers. In some sense, the algebra $\Zh$ might be thought of as a stably finite analogue of the Cuntz algebra $\Oh_{\infty}$ (cf.\ \cite{TW1} and \cite{TW2}). None of the above mentioned  counterexamples is $\Zh$-stable, but virtually all classes for which the Elliott conjecture has been confirmed so far consist of $\Zh$-stable $C^{*}$-algebras (\cite{TW2}). On the other hand, at the present stage it is not known how $\Zh$-stability alone can be used to obtain classification results, but in these notes we show that, in combination with other natural conditions, it can indeed be used to verify the Elliott conjecture for a large class of stably finite simple nuclear $C^{*}$-algebras.

$\Zh$-stable $C^{*}$-algebras behave very well in many respects (cf.\ \cite{GJS}, \cite{J} and \cite{R2}). For instance, a separable simple nuclear $\Zh$-stable $C^{*}$-algebra is either purely infinite or stably finite (in which case it has weakly unperforated $K_{0}$-group and satisfies Blackadar's second fundamental comparability property). Moreover, a stably finite separable simple $\Zh$-stable $C^{*}$-algebra has real rank zero if and only if the positive part of the ordered $K_{0}$-group has dense image in the positive continuous affine functions on the tracial state space. If one thinks of the tracial state space $T(A)$ as an underlying space of a $\Zh$-stable and stably finite nuclear $C^{*}$-algebra $A$, one might regard the elements of $A$ as continuous affine functions on $T(A)$. Then $A$ has real rank zero (i.e., elements with finite spectrum are dense in the subspace of all self-adjoint elements of $A$, cf.\ \cite{BP}) precisely if $K_{0}(A)_{+}$ is dense in $A_{+}$ in the (in general non-Hausdorff) topology coming from $T(A)$. In view of this observation it is natural to look for classification results in the real rank zero case first.

The known classification theorems for approximately homogeneous (AH) or approximately subhomogeneous (ASH) $C^{*}$-algebras all use conditions involving topological covering dimension of the underlying base spaces. This point of view was also pursued in \cite{W4}, where we verified the Elliott conjecture for separable simple unital $C^{*}$-algebras with finite decomposition rank and compact and zero-dimensional space of extremal tracial states. (The decomposition rank is a notion of covering dimension for nuclear $C^{*}$-algebras introduced in \cite{KW} by E.\ Kirchberg and the author; cf.\ also \cite{W1} and \cite{W2}.) Similar results were obtained in \cite{B} and \cite{L4}, under the additional assumption that the algebras have only countably many extremal tracial states. In the (locally) AH real rank zero case, the trace space conditions are not necessary, as was shown in \cite{L5} and \cite{Ng}. In the present article we show that the condition of \cite{W4} on the tracial state space can be removed if the algebras in question are $\Zh$-stable -- in other words, we verify the Elliott conjecture for the class of separable simple unital $\Zh$-stable $C^{*}$-algebras with real rank zero and finite decomposition rank (satisfying the UCT). Following the lines of \cite{W4}, we do not prove a classification result directly. Instead, we show that $\Zh$-stability, real rank zero and finite decomposition rank together imply tracial rank zero;  classification will then follow from a theorem of Lin.  As a consequence, our algebras have decomposition rank no greater than 2, and they are AH of topological dimension less than or equal to 3 and  ASH of topological dimension at most 2. Moreover, we conclude from \cite{EGL} and \cite{TW2} that, if $A$ has finite decomposition rank and real rank zero, then $A$ is $\Zh$-stable if and only if it is approximately divisible.

I would like to thank Ping Wong Ng, Andrew Toms and, in particular, Nate Brown for helpful comments and inspiring conversations.

\newpage

\section{Real rank zero, decomposition rank and order zero maps}

Below we recall the definitions of decomposition rank and of order zero maps and describe the particular situation for real rank zero algebras.

\bn
\label{d-dr} 
\begin{ndefn} (cf.\ \cite{KW}, Definitions 2.2 and 3.1) Let $A$ and $F$ be separable $C^*$-algebras, $F$ finite-dimensional.\\
(i) A completely positive map $\varphi :  F \to A$ has   order zero, $\ord \varphi = 0$, if it preserves orthogonality, i.e., $\varphi(e) \varphi(f) = \varphi(f) \varphi(e) = 0$ for all $e,f \in F$ with $ef = fe = 0$.\\ 
(ii) A completely positive map $\varphi : F \to A$ is $n$-decomposable, if there is a decomposition $F=F_{0} \oplus \ldots \oplus F_{n}$ such that the restriction of $\varphi$ to $F_{i}$ has   order zero for each $i \in \{0, \ldots, n\}$; we say $\varphi$ is $n$-decomposable with respect to $F=F_{0} \oplus \ldots \oplus F_{n}$.\\
(iii) $A$ has decomposition rank $n$, $\dr A = n$, if $n$ is the least integer such that the following holds: Given $\{b_1, \ldots, b_m\} \subset A$ and $\varepsilon > 0$, there is a completely positive approximation $(\bar{F}, \psi, \varphi)$ for $b_1, \ldots, b_m$ within $\varepsilon$ (i.e., $\bar{F}$ is a finite-dimensional $C^{*}$-algebra and $\psi:A \to \bar{F}$ and $\varphi:\bar{F} \to A$ are completely positive contractive and $\|\varphi \psi (b_i) - b_i\| < \varepsilon$) such that $\varphi$ is $n$-decomposable. If no such $n$ exists, we write $\dr A = \infty$.  
\end{ndefn}
\en

\bn
\label{order-zero}
Let $\varphi: F \to A$ be an   order zero map. We say a $*$-homomorphism $\sigma_{\varphi}:F \to A''$ is a supporting $*$-homomorphism for $\varphi$, if 
\begin{equation}
\label{central-support}
\varphi (x) = \varphi(\be_F) \sigma_{\varphi} (x) = \sigma_{\varphi}(x) \varphi(\be_F) \; \forall \, x \in F \, .
\end{equation}
By \cite{W4}, 1.2, any   order zero map has a supporting $*$-homomorphism. Moreover, each positive contractive $q \in C^*(\varphi(\be_F))$ defines a completely positive contractive   order zero map $\varphi_{q}: F \to A$ by 
\[
\varphi_{q}(\,.\,):= q \sigma_{\varphi}(\, .\,) \, .
\]
It is clear from (\ref{central-support}) that $q$ commutes with $\varphi(F)$. If $q$ happens to be a projection, then $\varphi_{q}$ is a $*$-homomorphism by \cite{W1}, Proposition 3.2. 
\en

\bn
\label{discrete-order-zero}
Following \cite{W4}, Definition 2.2, we say that a completely positive map 
\[
\varphi: F=M_{r_{1}} \oplus \ldots \oplus M_{r_{s}} \to A 
\] 
is a discrete order zero map, if $\ord \varphi = 0$ and each $\varphi(\be_{M_{r_{i}}})$, $i=1, \ldots, s$, is a multiple of a projection. The next lemma says that, if $A$ has real rank zero, then any order zero map into $A$ can be approximated by discrete order zero maps.

\begin{nlemma}
(\cite{W4}, Lemma 2.4) Let $A$ and $F$ be $C^{*}$-algebras, $A$ with real rank zero and $F$ finite-dimensional. Suppose $\varphi:F \to A$ is completely positive contractive with   order zero and let $\delta>0$ be given. Then there are a unital embedding $\iota:F \to \tilde{F}$ of $F$ into some finite-dimensional $C^{*}$-algebra $\tilde{F}$ and a discrete order zero map $\tilde{\varphi}:\tilde{F} \to A$ such that $\tilde{\varphi}(\be_{\tilde{F}})\le \varphi(\be_{F})$ and $\|\varphi(x) - \tilde{\varphi} \circ \iota(x)\|< \delta \cdot \|x\|$ for all $0 \neq x \in F$. 
\end{nlemma}
\en

\bn{\label{rr0dr}}
\begin{nprop} (\cite{W4}, Proposition 2.5) Let $A$ be a $C^{*}$-algebra with real rank zero and decomposition rank $n$. For any finite subset $\{b_{1},\ldots,b_{k}\} \in A$ and $\varepsilon>0$, there is a c.p.\ approximation $(F,\psi,\varphi)$ such that $\varphi$ is discretely $n$-decomposable and 
\[
\|\varphi\psi(b_{j})-b_{j}\|,\, \|\psi(b_{j})\psi(b_{l})-\psi({b_{j}b_{l}})\|<\varepsilon  \mbox{ for } j,\, l=1,\ldots,k \, .
\]
\end{nprop}
\en

\section{The Jiang--Su algebra and approximate tracial division}

Below we recall some facts about the Jiang--Su algebra $\Zh$ and establish one of the technical tools for the proof of our main theorem. 

\bn
\label{d-Z}
For relatively prime natural numbers $p$ and $q$ define their prime dimension drop interval to be the $C^{*}$-algebra  
\begin{displaymath}
I[p,q] := \{f \in \Ch([0,1],M_p \otimes M_{q}) \, | \, f(0) \in \be_{M_{p}} \otimes M_{q},  \,
f(1) \in M_{p} \otimes \be_{M_{q}} \} \, .
\end{displaymath}
In \cite{JS1}, Jiang and Su constructed an infinite-dimensional simple unital $C^{*}$-algebra $\Zh$, which  is $KK$-equivalent to $\C$ and has a unique normalized trace $\bar{\tau}$. It can be written as an inductive limit of prime dimension drop intervals; moreover, there is a unital embedding of any prime dimension drop interval into 
$\Zh$.
\en

\bn
\label{absorbing}
The Jiang--Su algebra is strongly self-absorbing in the sense of \cite{TW1}. In particular this means that, if $A$ is unital and $\Zh$-stable, there exists a 
sequence 
\[
(\theta_{n}: A \otimes \Zh \to A)_{n \in \N}
\]
of unital $*$-homomorphisms which  satisfies 
\[
\|\theta_{n}(a \otimes \be_{\Zh}) - a\| \stackrel{n \to \infty}{\longrightarrow} 0 \; \forall \, a \in A \, 
\] 
(cf.\ \cite{TW1}, Remark 2.7).
\en

\bn
\label{rr0>wu}
Let $A$ be simple and unital. By \cite{GJS}, $A$ and $A \otimes \Zh$ have isomorphic Elliott invariants iff the ordered group $K_{0}(A)$ is weakly unperforated. In \cite{T}, Toms has given an example of a simple unital $C^{*}$-algebra with weakly unperforated  $K_{0}$-group which is not $\Zh$-stable; this is a counterexample to the Elliott conjecture.  It is not known, although well possible, if real rank zero and finite decomposition rank together do imply $\Zh$-stability. At least we have the following

\begin{ntheorem}
Let $A$ be a separable simple unital $C^{*}$-algebra with real rank zero and finite decomposition rank. Then, $A$ and $A \otimes \Zh$ have isomorphic Elliott invariants; furthermore, $A \otimes \Zh$ has real rank zero and finite decomposition rank. 
\end{ntheorem}

\begin{nproof}
By \cite{W4}, Theorem 3.9, $K_{0}(A)$ is weakly unperforated, so $A$ and $A \otimes \Zh$ have isomorphic Elliott invariants. By \cite{W4}, Corollary 5.2, $K_{0}(A)_{+}$ is dense in $\Aff(T(A))_{+}$ (the positive continuous affine functions on the Choquet simplex $T(A)$), whence  $K_{0}(A \otimes \Zh)_{+}$ is dense in $\Aff(T(A \otimes \Zh))_{+}$. But now it follows from \cite{R2}, Theorem 7.2, that $A \otimes \Zh$ has real rank zero; $\dr (A \otimes \Zh)$ is finite by \cite{KW}, 3.2.
\end{nproof}
\en

\bn
\label{tracial-divisors}
The next two lemmas will provide a method of dividing elements of a $\Zh$-stable $C^{*}$-algebra approximately with respect to traces. They mark a key step in the proof of our main result, Theorem \ref{Z-stable>TAF}.

\begin{nlemma}
For any $n \in \N$ and $0< \beta <  1/2(n+1)$ there is a completely positive contractive   order zero map $\varrho: \C^{n+1} \to \Zh$ such that $\bar{\tau}(\varrho(e_{i})) > \beta$ for $i=1, \ldots, n+1$, where  the $e_{i}$ denote the canonical generators of $\C^{n+1}$. 
\end{nlemma}

\begin{nproof}
Choose $k \in \N$ such that 
\[
\frac{1}{k}< \frac{1}{2} \left(\frac{1}{2(n+1)} -\beta \right) \, .
\]
Set  $p:=k \cdot (n+1)$ and $q:=p+1$. Let $\kappa:\Ch([0,1]) \to I[p,q]$ be the canonical unital embedding and $\iota: I[p,q] \to \Zh$ a unital $*$-homomorphism; such an $\iota$ always exists (cf.\ \ref{d-Z}). 

Since $\bar{\tau} \circ \iota \circ \kappa$ is a state on $\Ch([0,1])$, it is straightforward to find orthogonal positive normalized functions $f_{0}$ and $f_{1}$ in $\Ch([0,1])$ such that 
\[
\bar{\tau} \circ \iota \circ \kappa(f_{0}+f_{1}) > \frac{1}{2} - \frac{n+1}{2}\left(\frac{1}{2(n+1)}-\beta \right)
\]
and $f_{0}(1)=f_{1}(0)=0$. 

Next, define a completely positive contractive map $\varphi: M_{q} \oplus M_{p} \to I[p,q]$ by 
\[
\varphi:= f_{0} \cdot (\be_{M_{p}} \otimes \id_{M_{q}}) + f_{1} \cdot (\id_{M_{p}} \otimes \be_{M_{q}}) \, ;
\]
it is straightforward to check that $\varphi$ has order zero. 

Let $\varrho_{0}: \C^{n+1} \to M_{q}$ and $\varrho_{1}: \C^{n+1} \to M_{p}$ be $*$-homomorphisms such that the projections $\varrho_{j}(e_{i})$ have rank $k$ for $j=0,1$ and $i=1, \ldots, n+1$. Define a completely positive contractive map $\varrho:\C^{n+1} \to I[p,q]$ by
\[
\varrho:= \varphi \circ (\varrho_{0} \oplus \varrho_{1}) \, .
\]
Being the composition of an order zero map and a $*$-homomorphism, $\varrho$ again has order zero.

If $g_{1}$, $g_{2} \in M_{q}$ and $h_{1}$, $h_{2} \in M_{p}$ are projections such that $\mathrm{rank} (g_{1}) = \mathrm{rank} (g_{2})$ and $\mathrm{rank}(h_{1}) = \mathrm{rank} (h_{2})$, there are partial isometries $v \in M_{q}$ and $s \in M_{p}$ such that $g_{1}=vv^{*}$, $g_{2} = v^{*}v$, $h_{1}=ss^{*}$ and $h_{2}=s^{*}s$. It is then straightforward to check that
\begin{eqnarray*}
\bar{\tau}\circ \varphi(g_{1},h_{1}) & = & \bar{\tau} \circ \varphi(g_{2},h_{2})\\
& = & \frac{\mathrm{rank}(g_{j})}{q} \cdot \bar{\tau} \circ \iota \circ \kappa(f_{0}) + \frac{\mathrm{rank}(h_{j})}{p} \cdot \bar{\tau} \circ \iota \circ \kappa(f_{1})
\, .
\end{eqnarray*}
As a consequence, we see that
\begin{eqnarray*}
\bar{\tau}\circ \varrho (e_{i}) & = & \frac{k}{k(n+1)+1} \cdot \bar{\tau} \circ \iota \circ \kappa(f_{0}) + \frac{1}{n+1} \cdot \bar{\tau} \circ \iota \circ \kappa(f_{1})\\
& \ge & \frac{1}{n+1} \cdot \bar{\tau} \circ \iota \circ \kappa(f_{0}) + \frac{1}{n+1} \cdot \bar{\tau} \circ \iota \circ \kappa(f_{1})  - \frac{1}{k} \\
& > & \frac{1}{2(n+1)} -  \frac{1}{2}\left(\frac{1}{2(n+1)}-\beta \right)- \frac{1}{k} \\
& > & \beta
\end{eqnarray*}
for $i=1, \ldots, n+1$.
\end{nproof}
\en

\bn
\label{order0-truncation}
\begin{nlemma}
Let $\varphi: F \to A$ and $\varrho: \C^{n+1} \to B$ be  completely positive contractive maps; suppose $\varphi$ is $n$-decomposable with respect to the decomposition $F= \bigoplus_{i=0}^{n} F_{i}$ and $\varrho$ has   order zero.\\
Then the completely positive contractive  map $\bar{\varphi}: F \to A \otimes B$, given by
\[
\bar{\varphi} (x):= \sum_{i=0}^{n} \varphi(x \be_{F_{i}}) \otimes \varrho(e_{i+1}) \; \forall \, x \in F \, ,
\]
has   order zero.
\end{nlemma}

\begin{nproof}
Obvious.
\end{nproof}
\en

\section{Tracial rank zero}

In this section we recall the notion of tracially AF $C^{*}$-algebras and provide an alternative characterization of when a simple unital $C^{*}$-algebra with Blackadar's second fundamental comparability property is tracially AF.

\bn
\label{d-tr0}
In \cite{L1}, Lin introduced the notion of tracially AF algebras (or, equivalently, $C^{*}$-algebras of tracial rank zero). In the simple and unital case, the definition reads as follows:
 
\begin{ndefn}
A simple unital $C^{*}$-algebra $A$ is said to be tracially AF, if, for any finite subset $\Fh \subset A$, $\varepsilon > 0$ and $0 \neq a \in A_{+}$, there is a finite-dimensional $C^{*}$-subalgebra $B \subset A$ with the following properties:
\begin{itemize}
\item[(i)] $\|b \be_{B} - \be_{B} b\|<\varepsilon \; \forall \, b \in \Fh$ 
\item[(ii)] $\dist(\be_{B}b\be_{B}, B)< \varepsilon \; \forall \, b \in \Fh$ 
\item[(iii)] $\be_{A}-\be_{B}$ is Murray--von Neumann equivalent to a projection in the hereditary subalgebra $\overline{aAa}$ of $A$.
\end{itemize}
\end{ndefn}
\en

\bn
\label{characterization-tr0}
Recall that a unital $C^{*}$-algebra $A$ is said to have Blackadar's second fundamental comparability property, if,  whenever $p$ and $q$ are projections in $A$ such that  $\tau(p) < \tau(q) \; \forall \, \tau \in T(A)$, then $p$ is Murray--von Neumann equivalent to a subprojection of $q$. In \cite{Li3}, Corollary 6.15, Lin provided a characterization of when a $C^{*}$-algebra with real rank zero and comparability is tracially AF. For our purposes the following slightly different version of Lin's characterization will be more useful. It is probably well known; but since we couldn't find an explicit  proof in the literature, we include one here.  

\begin{nlemma}
Let $A$ be a simple unital $C^{*}$-algebra which has the second fundamental comparability property; suppose that $K_{0}(A)_{+}$ has dense image in the positive continuous affine functions on $T(A)$, $\Aff(T(A))_{+}$, and that every hereditary subalgebra of $A$ contains a nonzero projection. Then $A$ is tracially AF if and only if there is $n \in \N$ such that, for any finite subset $\Fh \subset A$ and $\varepsilon > 0$, there is a finite-dimensional $C^{*}$-subalgebra $B \subset A$ with the following properties:
\begin{itemize}
\item[(i)] $\|b \be_{B} - \be_{B} b\|<\varepsilon \; \forall \, b \in \Fh$ 
\item[(ii)] $\dist(\be_{B}b\be_{B}, B)< \varepsilon \; \forall \, b \in \Fh$ 
\item[(iii)] $\tau(\be_{B})>\frac{1}{n} \; \forall \, \tau \in T(A)$.
\end{itemize}
\end{nlemma}

As for the hypotheses of the preceding lemma, note that if a $C^{*}$-algebra $A$ has comparability and real rank zero, then automatically $K_{0}(A)_{+}$ has dense image in $\Aff(T(A))_{+}$ and every hereditary subalgebra of $A$ contains a nonzero projection.\\
Before we turn to the proof of the lemma, we need two more technical results.
\en

\bn
\label{weakly-stable}
Recall that any finite-dimensional $C^{*}$-algebra can be written as the universal $C^{*}$-algebra generated by a set of matrix units with the respective relations. By the results of Chapter 14 in \cite{Lo}, finite-dimensional $C^{*}$-algebras are semiprojective. Equivalently, the defining relations are stable (cf.\ \cite{Lo}, Definition 14.1.1 and Theorem 14.1.4). As a direct consequence we obtain the following:
  
\begin{nprop}
Let $B$ be a finite-dimensional $C^{*}$-algebra. For every $\gamma>0$ there is $\vartheta>0$ such that the following holds: \\
Suppose $A$ is another $C^{*}$-algebra and $\varphi: B \to A$ is a $*$-homomorphism. If $p \in A$ is a projection satisfying 
\[
\|p \varphi(b) - \varphi(b) p\|< \vartheta \cdot \|b\| \; \forall \, 0 \neq b \in B \, ,
\]
then there is a $*$-homomorphism $\tilde{\varphi}: B \to pAp$ such that 
\[
\|\tilde{\varphi}(b) - p\varphi(b)p \| < \gamma \cdot \|b\| \; \forall \, 0 \neq b \in B \, .
\]
\end{nprop}
\en

\bn
\label{weakly-stable-2}
\begin{nlemma}
Let $A$ be a simple unital $C^{*}$-algebra which has the second fundamental comparability property; suppose that $K_{0}(A)_{+}$ has dense image in the positive continuous affine functions on $T(A)$. Let $\alpha, \beta >0$ and  a nontrivial  projection $p \in A$ with $\tau(p)>\alpha + 2 \beta \; \forall \, \tau \in T(A)$ be given.\\ 
Then, there are $s, t \in \N$ and $*$-homomorphisms $\iota_{0}:M_{t} \oplus M_{t+1} \to M_{t+s} \oplus M_{t+1+s}$, $\iota_{1}: M_{s} \oplus M_{s} \to M_{t+s} \oplus M_{t+1+s}$ and $\kappa: M_{t+s} \oplus M_{t+1+s} \to A$ such that the map $\iota_{0} + \iota_{1}$ is a unital $*$-homomorphism, 
\[
\kappa \iota_{0}(\be_{M_{t}}+\be_{M_{t+1}}) = \be_{A} -p \, ,
\]
\[
\frac{s}{t+1} > \frac{\alpha}{1-\alpha} + \beta 
\]
and
\[
t \cdot \theta \iota_{1}(\be_{M_{s} \oplus M_{s}}) < s \cdot \theta \iota_{0}(\be_{M_{t} \oplus M_{t+1}}) < (t+1) \cdot \theta \iota_{1} (\be_{M_{s} \oplus M_{s}})
\]
for any tracial state $\theta$ on $M_{t+s} \oplus M_{t+1+s}$.
\end{nlemma}

\begin{nproof}
Let $f \in \Aff (T(A))_{+}$ be the image of $\be_{A} - p$ under the canonical map $r: K_{0}(A) \to \Aff (T(A))$. Note that $f$ is nowhere zero since $A$ is simple and any tracial state on $A$ is faithful.

We first show that, for any $t \in \N$, there is a unital $*$-homomorphism 
\[
\nu : M_{t} \oplus M_{t+1} \to (\be_{A}-p)A(\be_{A}-p) \, .
\]
By our hypotheses on $A$, there is a projection $e \in (\be_{A}-p)A(\be_{A}-p)$ such that 
\[
\frac{1}{t+1} \cdot f < r(e) < \frac{1}{t} \cdot f \, .
\]
From comparison, we now obtain $t$ pairwise orthogonal subprojections of $\be_{A}-p$ which are all Murray--von Neumann equivalent to $e$; this yields a $*$-homomorphism 
\[
\bar{\nu}: M_{t} \to (\be_{A}-p)A(\be_{A}-p) \, .
\] 
We have
\[
r(\be_{A}-p-\bar{\nu}(\be_{M_{t}})) = f - t \cdot r(e) < r(e) \, ,
\]
whence there is a partial isometry $s \in A$ such that 
\[
s^{*}s \le \bar{\nu}(e_{11}) 
\]
and
\[
s s^{*} = \be_{A} - p - \bar{\nu}(\be_{M_{t}}) \, .
\]
But now it is straightforward to construct $*$-homomorphisms 
\[
\nu_{1}: M_{t+1} \to (\be_{A}-p)A(\be_{A}-p)
\]
and 
\[
\nu_{0}: M_{t} \to (\be_{A}-p)A(\be_{A}-p)
\]
such that 
\[
\nu_{1}(e_{t+1,t+1}) = \be_{A} - p - \bar{\nu}(\be_{M_{t}})
\]
and
\[
\nu_{0}(e_{jj}) + \nu_{1}(e_{jj}) = \bar{\nu}(e_{jj}) \; \forall \, j \in \{1, \ldots , t\} \, .
\]
The map $\nu_{1}$ is determined by setting $\nu_{1}(e_{t+1,1}):= (\be_{A} - p - \bar{\nu}(\be_{M_{t}})) s$ and $\nu_{1}(e_{j,1}) := \bar{\nu}(e_{j,1}) s^{*}s$ for $j=1, \ldots,t$.  $\nu_{0}$ is determined by $\nu_{0}(e_{j,1}):= \bar{\nu}(e_{j,1}) (\bar{\nu}(e_{11})-s^{*}s)$ for $j=1, \ldots,t$.  The maps $\nu_{0}$ and $\nu_{1}$ have orthogonal images and therefore add up to a unital $*$-homomorphism
\[
\nu: M_{t} \oplus M_{t+1} \to (\be_{A}-p)A(\be_{A}-p) \, .
\]

Next, one observes that 
\[
\frac{\alpha + 2 \beta}{1 - \alpha - 2 \beta} > \frac{\alpha + \beta}{1-\alpha-\beta}> \frac{\alpha}{1-\alpha} + \beta \, ;
\]
we may therefore fix $s,t \in \N$ such that 
\[
\frac{\alpha + 2 \beta}{1 - \alpha - 2 \beta} > \frac{s}{t} > \frac{s}{t+1} > \frac{\alpha}{1-\alpha} + \beta \, .
\]
By hypothesis, we have
\[
\tau(\be_{A}-p) < 1 - \alpha - 2 \beta \; \forall \, \tau \in T(A) \, ,
\]
so
\[
\frac{s}{t} \cdot \tau(\be_{A}-p) < \alpha + 2 \beta < \tau(p) \; \forall \, \tau \in T(A) \, .
\]
Since
\[
r(\nu_{0}(e_{11}) + \nu_{1}(e_{11})) < \frac{1}{t} \cdot f = \frac{1}{t} \cdot r(\be_{A} - p) \, ,
\]
from comparison we obtain $s$ pairwise orthogonal subprojections of $p$ which are all Murray--von Neumann equivalent to the projection $\nu_{0}(e_{11})+ \nu_{1}(e_{11})$. This clearly yields an embedding
\[
\kappa : M_{t+s} \oplus M_{t+1+s} \to A
\]
with the desired properties, where $\iota_{0}$ is the sum of the upper left corner embeddings and $\iota_{1}$ is the sum of the lower right corner embeddings.
\end{nproof}
\en

\begin{nproof}
(of Lemma \ref{characterization-tr0}) If $A$ is tracially AF, the assertion obviously holds with any $n>1$. Conversely, suppose there is $n \in \N$ such that the assertion in the lemma holds. Let $0\neq a \in A_{+}$, a finite subset $\Fh \subset A$ and $\varepsilon>0$ be given. By our assumption on $A$, there is a nonzero projection $q \in \overline{aAa}$. Set 
\[
\eta:= \min\{\tau(q) \, | \, \tau \in T(A)\} \, .
\]
Note that $\eta>0$, since $T(A)$ is compact and the positive function $\tau \mapsto \tau(q)$ is continuous and everywhere nonzero ($A$ is simple, hence every tracial state is faithful).

For $i \in \N$, choose strictly positive numbers $\varepsilon_{i}$ such that $\sum_{\N} \varepsilon_{i} < \varepsilon$ and define 
\[
\alpha_{i}:= \frac{1}{n}  \sum_{k=0}^{i} \left(1 - \frac{1}{n} \right)^{k} \, ;
\]
note that $(1- \alpha_{i})(1-\frac{1}{n}) = 1- \alpha_{i+1}$, whence 
\begin{equation}
\label{alpha_i}
\frac{\alpha_{i}}{1-\alpha_{i}} (1 - \alpha_{i+1}) + (1 - \alpha_{i+1}) = 1- \frac{1}{n} \; \forall \, i
\end{equation}
and that 
\[
\alpha_{i} \stackrel{i \to \infty}{\longrightarrow} 1 \, .
\]
We shall inductively construct finite-dimensional $C^{*}$-algebras $B_{i} \subset A$, $i \in \N$, with the following properties:
\begin{itemize}
\item[a)] $\|b \be_{B_{i}} - \be_{B_{i}} b\| < \sum_{k=0}^{i} \varepsilon_{k} \; \forall \, b \in \Fh$
\item[b)] $\dist(\be_{B_{i}}b \be_{B_{i}}, B_{i}) < \sum_{k=0}^{i} \varepsilon_{k} \; \forall \, b \in \Fh$
\item[c)] $\tau(\be_{B_{i}}) > \alpha_{i} \; \forall \, \tau \in T(A)$.
\end{itemize}
Having done so, (i) and (ii) of Definition \ref{d-tr0} will hold by a) and b), respectively. Since $\alpha_{i} \to 1$, by c) there is $K \in \N$ such that $\tau(\be_{A} - \be_{B_{K}}) < \eta$. From comparison we see that $\be_{A} - \be_{B_{K}}$ is Murray--von Neumann equivalent to  $q \in \overline{aAa}$, whence (iii) of \ref{d-tr0} holds and $A$ is tracially AF. 

We proceed to construct the $B_{i}$. $B_{0}$ obviously exists by assumption. Suppose that $B_{i}$ (satisfying properties a), b) and c)) has been constructed for some $i \in \N$. We show how to obtain $B_{i+1}$.\\
First, set
\[
\beta := \halb \min \{ \tau(\be_{B_{i}}) - \alpha_{i} \, | \, \tau \in T(A) \} \, ,
\]
then $\beta >0$, again since $T(A)$ is compact and the positive function $\tau \mapsto \tau(\be_{B_{i}})- \alpha_{i}$ is continuous and nonzero. \\
By Lemma \ref{weakly-stable-2}, there are $s, t \in \N$, finite-dimensional $C^{*}$-algebras $D:=M_{t+s} \oplus M_{t+1+s}$, $C_{0}:= M_{t} \oplus M_{t+1}$ and $C_{1}:= M_{s} \oplus M_{s}$ and $*$-homomorphisms $\iota_{0}:C_{0} \to D$, $\iota_{1}: C_{1} \to D$ and $\kappa: D \to A$ satisfying
\begin{equation}
\label{s-t+1}
\frac{s}{t+1} > \frac{\alpha_{i}}{1-\alpha_{i}} + \beta \, ,
\end{equation}
\begin{equation}
\label{unit-iota_0}
\kappa \iota_{0}(\be_{C_{0}}) = \be_{A} - \be_{B_{i}} \, ,
\end{equation}
\begin{equation*}
\iota_{0}(\be_{C_{0}}) + \iota_{1}(\be_{C_{1}}) = \be_{D}
\end{equation*} 
and 
\begin{equation}
\label{trace-ratio}
t \cdot \theta \iota_{1}(\be_{C_{1}}) < s \cdot \theta \iota_{0}(\be_{C_{0}}) < (t+1) \cdot \theta \iota_{1}(\be_{C_{1}})
\end{equation}
for any $\theta \in T(D)$. In particular,
\begin{equation}
\label{kappaiota}
\kappa \iota_{1}(\be_{C_{1}}) \le \be_{B_{i}} \, .
\end{equation}
Set 
\[
\widetilde{\Fh} := \Fh \cup \Bh_{1}(B_{i}) \cup \Bh_{1}(\kappa(D)) 
\]
(where $\Bh_{1}( \, . \, )$ denotes the unit ball of a $C^{*}$-algebra); choose $\gamma>0$ such that 
\begin{equation}
\label{alphabeta}
2 \gamma \left(\frac{\alpha_{i}}{1- \alpha_{i}} + \beta +1 \right) < \beta (1 - \alpha_{i+1})
\end{equation}
and
\begin{equation}
\label{gammaepsilon}
7 \gamma < \varepsilon_{i+1} \, .
\end{equation}
Choose $0 < \vartheta < \gamma$ such that the assertion of Proposition \ref{weakly-stable} holds for both the finite-dimensional $C^{*}$-algebras $B_{i}$ and $D$.

By hypothesis, there is a finite-dimensional $C^{*}$-algebra $F \subset A$ such that 
\begin{itemize}
\item[d)] $\|b \be_{F} - \be_{F} b \| < \vartheta \; \forall \, b \in \widetilde{\Fh}$
\item[e)] $\dist( \be_{F} b \be_{F}, F) < \vartheta \; \forall \, b \in \widetilde{F}$
\item[f)] $\tau(\be_{F}) > \frac{1}{n} \; \forall \, \tau \in T(A)$.
\end{itemize}
By d) in connection with Proposition \ref{weakly-stable}, there is a $*$-homomorphism 
\[
\varrho: B_{i} \to (\be_{A}-\be_{F}) A (\be_{A} - \be_{F})
\]
such that 
\begin{equation}
\label{perturbed-varrho}
\|\varrho(b) -  (\be_{A}-\be_{F}) b (\be_{A} - \be_{F}) \| < \gamma \cdot \|b\| \; \forall \, 0 \neq b \in B_{i} \, .
\end{equation}
Similarly, we obtain a $*$-homomorphism 
\[
\tilde{\kappa}: D \to (\be_{A} - \be_{F}) A (\be_{A} - \be_{F})
\]
such that
\begin{equation}
\label{perturbed-kappa}
\|\tilde{\kappa}(d) -  (\be_{A}-\be_{F}) \kappa(d) (\be_{A} - \be_{F}) \| < \gamma \cdot \|d\| \; \forall \, 0 \neq d \in D \, .
\end{equation}
Set 
\[
B_{i+1}:= \varrho(B_{i}) \oplus F \, ;
\]
we proceed to check that $B_{i+1}$ statisfies properties a), b) and c) above (with $i+1$ in place of $i$).

First, we have
\begin{eqnarray*}
\lefteqn{\|b \be_{B_{i+1}} - \be_{B_{i+1}} b \|}\\
& = & \|b (\be_{F} + \varrho(\be_{B_{i}})) - (\be_{F} + \varrho(\be_{B_{i}})) b \| \\
& \le & \|b \varrho(\be_{B_{i}}) - \varrho(\be_{B_{i}}) b \| + \|b \be_{F} - \be_{F} b \| \\
& \stackrel{\mathrm{d)}, (\ref{perturbed-varrho})}{<} & \|b (\be_{A} - \be_{F}) \be_{B_{i}} (\be_{A} - \be_{F}) - (\be_{A} - \be_{F}) \be_{B_{i}} (\be_{A} - \be_{F}) b \| \\
& &  + 2 \gamma + \vartheta\\
& \stackrel{\mathrm{d)}}{<} & \| (\be_{A} - \be_{F}) (b \be_{B_{i}} - \be_{B_{i}} b) (\be_{A} - \be_{F}) \| \\
& &  + 2 \gamma + 3 \vartheta\\
& \stackrel{(\mathrm{a}),(\ref{gammaepsilon})}< & \sum_{k=0}^{i+1} \varepsilon_{k} \, ,
\end{eqnarray*}
so a) above holds  (with $i+1$ in place of $i$). Next, we check b):
\begin{eqnarray*}
\lefteqn{\dist( \be_{B_{i+1}} b \be_{B_{i+1}}, B_{i+1})} \\
& = & \dist((\varrho(\be_{B_{i}}) + \be_{F}) b (\varrho(\be_{B_{i}}) + \be_{F}), B_{i+1} ) \\
& \stackrel{\mathrm{d)}}{\le} & \dist(\varrho(\be_{B_{i}}) b \varrho(\be_{B_{i}}) + \be_{F} b \be_{F}, B_{i+1} ) + 2 \vartheta \\
& \le & \dist(\varrho(\be_{B_{i}}) b \varrho(\be_{B_{i}}) , \varrho(B_{i}))  + \dist(\be_{F}b \be_{F}, F)  + 2 \vartheta \\
& = & \dist(\varrho(\be_{B_{i}}) (\be_{A}-\be_{F}) b (\be_{A}-\be_{F}) \varrho(\be_{B_{i}}) , \varrho(B_{i}))  + \dist(\be_{F}b \be_{F}, F)  + 2 \vartheta \\
& \stackrel{(\ref{perturbed-varrho})}{<} & \dist(\varrho(\be_{B_{i}}) \varrho( b) \varrho(\be_{B_{i}}) , \varrho(B_{i}))  + \dist(\be_{F}b \be_{F}, F)  + \gamma + 2 \vartheta \\
& \stackrel{\mathrm{b)},\mathrm{e)}}{<} & \sum_{k=0}^{i} \varepsilon_{k} +  \gamma + 3 \vartheta \\
& \stackrel{(\ref{gammaepsilon})}{<} & \sum_{k=0}^{i+1} \varepsilon_{k} \, .
\end{eqnarray*}

We will not prove c) directly; instead, we assume that
\begin{equation}
\label{tau-assumption}
\tau(\be_{B_{i+1}}) \le \alpha_{i+1}
\end{equation}
for some $\tau \in T(A)$ to deduce a contradiction. 

\noindent
We then have
\begin{eqnarray}
\lefteqn{\tau((\be_{A} - \be_{F}) \be_{B_{i}}(\be_{A} - \be_{F}))} \nonumber \\
& \stackrel{(\ref{kappaiota})}{\ge} & \tau((\be_{A} - \be_{F}) \kappa \iota_{1}(\be_{C_{1}}) (\be_{A} - \be_{F})) \nonumber \\
& \stackrel{(\ref{perturbed-kappa})}{>} & \tau \tilde{\kappa} \iota_{1}(\be_{C_{1}}) - \gamma \nonumber \\
& \stackrel{(\ref{trace-ratio})}{>} & \frac{s}{t+1} \cdot \tau \tilde{\kappa} \iota_{0}(\be_{C_{0}}) - \gamma \nonumber \\
& \stackrel{(\ref{perturbed-kappa})}{>} & \frac{s}{t+1} \cdot (\tau((\be_{A}-\be_{F}) \kappa \iota_{0}(\be_{C_{0}}) (\be_{A}-\be_{F})) - \gamma) - \gamma \nonumber \\
& \stackrel{(\ref{unit-iota_0})}{=} & \frac{s}{t+1} \cdot (\tau((\be_{A}-\be_{F}) (\be_{A} - \be_{B_{i}}) (\be_{A}-\be_{F})) - \gamma) - \gamma \nonumber \\
& \stackrel{(\ref{perturbed-varrho})}{>} & \frac{s}{t+1} \cdot (\tau(\be_{A} - \be_{F} - \varrho(\be_{B_{i}})) - 2 \gamma) - \gamma \nonumber \\
& = & \frac{s}{t+1} \cdot (\tau(\be_{A} - \be_{B_{i+1}}) - 2 \gamma) - \gamma  \nonumber \\
& \stackrel{(\ref{tau-assumption})}{\ge} & \frac{s}{t+1} \cdot (1 - \alpha_{i+1} - 2 \gamma) - \gamma \nonumber \\
& \stackrel{(\ref{s-t+1})}{>} & \left(\frac{\alpha_{i}}{1 - \alpha_{i}} + \beta \right) (1 - \alpha_{i+1} - 2 \gamma) - \gamma \, . \label{trace-estimate}
\end{eqnarray}
As a consequence, we obtain
\begin{eqnarray*}
\lefteqn{\tau(\be_{A} - \be_{F})}\\
& = & \tau((\be_{A}- \be_{F}) \be_{B_{i}}(\be_{A}- \be_{F})) \\
& & + \tau((\be_{A}- \be_{F}) (\be_{A}-\be_{B_{i}})(\be_{A}- \be_{F})) \\
& \stackrel{(\ref{trace-estimate}),(\ref{perturbed-varrho})}{>} & \left(\frac{\alpha_{i}}{1-\alpha_{i}} + \beta \right) (1-\alpha_{i+1} - 2 \gamma) - \gamma \\
& & + \tau(\be_{A} - \be_{F} - \varrho(\be_{B_{i}})) - \gamma \\
& = & \left(\frac{\alpha_{i}}{1-\alpha_{i}} + \beta \right) (1-\alpha_{i+1} - 2 \gamma) - \gamma \\
& & + \tau(\be_{A} - \be_{B_{i+1}}) - \gamma \\
& \stackrel{(\ref{tau-assumption})}{\ge} & \left(\frac{\alpha_{i}}{1-\alpha_{i}} + \beta \right) (1-\alpha_{i+1} - 2 \gamma) - 2\gamma  + 1 - \alpha_{i+1} \\
& = & \frac{\alpha_{i}}{1-\alpha_{i}} (1 - \alpha_{i+1}) + (1 - \alpha_{i+1}) \\
& & - 2 \gamma  \left(\frac{\alpha_{i}}{1-\alpha_{i}} + \beta + 1 \right) + \beta (1 - \alpha_{i+1}) \\
& \stackrel{(\ref{alpha_i}),(\ref{alphabeta})}{>} & 1- \frac{1}{n} \, ,
\end{eqnarray*}
a contradiction to f). Therefore, (\ref{tau-assumption}) is wrong and we have
\[
\tau(\be_{B_{i+1}}) > \alpha_{i+1} \, ,
\]
whence c) holds.

Induction now yields finite-dimensional subalgebras $B_{i} \subset A$, $i \in \N$, with properties a), b) and c). We are done.
\end{nproof}

\section{The main result}

\bn
\label{Z-stable>TAF}
\begin{ntheorem}
Let $A$ be a separable simple unital $C^{*}$-algebra with finite decomposition rank $n$; suppose  $A$ is $\Zh$-stable and has real rank zero. Then, $A$ is tracially AF.
\end{ntheorem}
\en

\bn
\label{multiplicative-domain}
Before turning to the proof, we note the following consequence of Stinespring's theorem, which  is a standard tool to analyze completely positive approximations of nuclear $C^{*}$-algebras. See \cite{KW}, Lemma 3.5, for a proof.
 
\begin{nlemma}
Let $A$ and $F$ be $C^{*}$-algebras, $b \in A$ a normalized positive element and $\eta>0$. If $A \stackrel{\psi}{\longrightarrow} F \stackrel{\varphi}{\longrightarrow} A$ are completely positive contractive maps satisfying 
\[
\|\varphi \psi(b) - b\|, \, \|\varphi \psi(b^{2}) - b^{2}\| < \eta \, ,
\]
then, for any $0 \neq x \in F$, 
\[
\|\varphi(\psi(b)x)- \varphi \psi(b) \varphi(x)\|< 2 \eta^{\halb} \|x\| \, .
\]
\end{nlemma}
\en

\begin{nproof}
(of Theorem \ref{Z-stable>TAF})  We are going to show that $A$ satisfies the conditions of Lemma \ref{characterization-tr0} with $5(n +1)$ in place of $n$. For convenience, we define
\begin{equation}
\label{mu}
\mu:= \frac{1}{5(n+1)} \, .
\end{equation}

Let $\varepsilon>0$ and a finite subset $\Fh \subset A$  be given. We may assume that the elements of $\Fh$ are positive and normalized and that $\be_{A} \in \Fh$.  
Choose strictly positive numbers $\alpha$,  $\delta$  and  $\eta$ such that 
\begin{equation}
\label{constants}
\alpha + \eta + \delta< \mu \mbox{ and } 2 \alpha + 4 \frac{\delta+ \alpha^{\halb}}{\eta} < \frac{\varepsilon}{3} \, .
\end{equation}

Since $\dr A=n$, there is an $n$-decomposable completely positive contractive approximation $(F,\psi,\varphi)$  for $\Fh^{2}$ within $\alpha$. More precisely, there are a finite-dimensional $C^{*}$-algebra $F$ and completely positive contractive maps 
\[
A  \stackrel{\psi}{\longrightarrow} F \stackrel{\varphi}{\longrightarrow} A
\]
such that 
\begin{equation}
\label{cp-appr}
\|\varphi \psi (b_{1}b_{2}) - b_{1}b_{2}\|< \alpha \; \forall \, b_{1},b_{2} \in \Fh
\end{equation}
and such that $\varphi$ is $n$-decomposable with respect to the decomposition $F=F_{0} \oplus \ldots \oplus F_{n}$. By Proposition \ref{rr0dr} we may even assume that $\|\psi(b_{1}b_{2})- \psi(b_{1})\psi(b_{2})\|< \alpha$ for all $b_{1}, b_{2} \in \Fh$.

By Lemma \ref{tracial-divisors} there is a completely positive contractive   order zero map $\varrho: \C^{n+1} \to \Zh$ such that 
\begin{equation}
\label{bartaubeta}
\beta < \bar{\tau}(\varrho(e_{i}))
\end{equation} 
for $i=1, \ldots, n+1$, where 
\begin{equation}
\label{betamu}
\beta:= 2 \mu \, , 
\end{equation}
the $e_{i}$ again denote the canonical generators of $\C^{n+1}$ and $\bar{\tau}$ is the unique tracial state on $\Zh$. As in 
Lemma \ref{order0-truncation}, define a completely positive contractive    order zero map $\bar{\varphi}: F \to A \otimes \Zh$ by
\begin{equation}
\label{barphi}
\bar{\varphi} (x):= \sum_{i=0}^{n} \varphi(x \be_{F_{i}}) \otimes \varrho(e_{i+1}) \; \forall \, 0 \neq x \in F \, .
\end{equation}

From Lemma \ref{discrete-order-zero} we obtain a unital embedding $\iota:F \to \tilde{F}$ of $F$ into some finite-dimensional $C^{*}$-algebra $\tilde{F}$ and a discrete order zero map $\tilde{\varphi}:\tilde{F} \to A \otimes \Zh$ such that $\tilde{\varphi}(\be_{\tilde{F}})\le \bar{\varphi}(\be_{F})$ and 
\begin{equation}
\label{barphitildephi}
\|\bar{\varphi}(x) - \tilde{\varphi} \circ \iota(x)\|< \delta \cdot \|x\| \; \forall \, 0 \neq x \in F \, .
\end{equation} 

Let  $\chi_{\eta}$ denote the characteristic function on the interval $(\eta, \infty)$. The spectrum of $\tilde{\varphi}(\be_{\tilde{F}})$ is finite because $\tilde{\varphi}$ has discrete order zero, so  
\begin{equation}
\label{qdef}
q:= \chi_{\eta}(\tilde{\varphi}(\be_{\tilde{F}})) 
\end{equation}
defines a projection in $C^{*}(\tilde{\varphi}(\be_{\tilde{F}})) \subset A \otimes \Zh$. Let $\sigma_{\tilde{\varphi}}: \tilde{F} \to (A \otimes \Zh)''$ be a supporting $*$-homomorphism for $\tilde{\varphi}$.  As in \ref{order-zero} we can define a $*$-homomorphism $\sigma (= \tilde{\varphi}_{q}): \tilde{F} \to A  \otimes \Zh$ by 
\begin{equation}
\label{sigmadef}
\sigma(x):= q \sigma_{\tilde{\varphi}}(x) \; \forall \, x \in \tilde{F} \, ;
\end{equation}
we have
\begin{equation}
\label{qsigma}
\sigma(\be_{\tilde{F}}) = q \, .
\end{equation}
By functional calculus there is a positive element $h \in C^{*}(\tilde{\varphi}(\be_{\tilde{F}}))$ with 
\begin{equation}
\label{normh}
\|h\|\le \frac{1}{\eta}
\end{equation}
satisfying 
\begin{equation}
\label{qh}
q = q  h  \tilde{\varphi}(\be_{\tilde{F}}) \, . 
\end{equation}
Since $\tilde{\varphi}(x)= \tilde{\varphi}(\be_{\tilde{F}}) \sigma_{\tilde{\varphi}}(x)$, we have
\begin{eqnarray}
\label{sigma-factorization}
\sigma(x) = q h \tilde{\varphi}(x) \, ;
\end{eqnarray}
moreover,  
\begin{equation}
\label{qhtildephi}
[q,\tilde{\varphi}(\tilde{F})] = [h,\tilde{\varphi}(\tilde{F})] = [q,h] = 0
\end{equation} by \ref{order-zero}.

By \ref{absorbing} there is a unital $*$-homomorphism $\theta: A \otimes \Zh \to A$  satisfying
\begin{equation}
\label{btheta}
\|\theta(b \otimes \be_{\Zh}) - b\| < \frac{\varepsilon}{3}  \; \forall \, b \in \Fh  \, .
\end{equation}
We now define a finite-dimensional $C^{*}$-subalgebra of $A$ by
\begin{equation}
\label{Bdef}
B:= \theta \circ \sigma(\tilde{F}) \subset A  
\end{equation}
and proceed to check that it satisfies properties (i), (ii) and (iii) of Lemma \ref{characterization-tr0}.\\

For $b \in \Fh$ we have
\begin{eqnarray*}
\lefteqn{\|[\be_{B}\, , \, b] \|}\\
& \stackrel{(\ref{btheta}),(\ref{Bdef})}{\le} & \|\theta([\sigma(\be_{\tilde{F}}) \, , \, b \otimes \be_{\Zh}])\| + \frac{2 \varepsilon}{3}\\
& \le & \| [\sigma(\be_{\tilde{F}}) \, , \, b \otimes \be_{\Zh}] \| + \frac{2 \varepsilon}{3}\\
& \stackrel{(\ref{cp-appr})}{\le} & \| [\sigma(\be_{\tilde{F}}) \, , \, \varphi \psi(b) \otimes \be_{\Zh}] \| + \frac{2 \varepsilon}{3} + 2 \alpha\\
& \stackrel{(\ref{sigma-factorization})}{=} & \| [q h \tilde{\varphi}(\be_{\tilde{F}}) \, , \, \varphi \psi(b) \otimes \be_{\Zh}] \| + \frac{2 \varepsilon}{3} + 2 \alpha\\
& \stackrel{(\ref{normh}),(\ref{barphitildephi})}{\le} & \| q h \bar{\varphi}(\be_{F})  (\varphi \psi(b) \otimes \be_{\Zh}) -  (\varphi \psi(b) \otimes \be_{\Zh}) \bar{\varphi}(\be_{F}) q h\| + \frac{2 \varepsilon}{3} + 2 \alpha + \frac{2 \delta}{\eta}\\
& \stackrel{(\ref{barphi})}{=} & \| [q h , \sum_{i=0}^{n}( (\varphi (\be_{F_{i}})  \varphi \psi(b)) \otimes \varrho(e_{i+1}))] \|  + \frac{2 \varepsilon}{3} + 2 \alpha + \frac{2 \delta}{\eta}\\
& \stackrel{\ref{multiplicative-domain},(\ref{cp-appr})}{\le} & \| [q h, \sum_{i=0}^{n} (\varphi (\be_{F_{i}}   \psi(b)) \otimes \varrho(e_{i+1}))\|  + \frac{2 \varepsilon}{3} + 2 \alpha + \frac{2 \delta}{\eta} + \frac{4 \alpha^{\halb}}{\eta}\\
& \stackrel{(\ref{barphi})}{=} & \|[q h \, , \, \bar{\varphi}\psi(b)] \|  + \frac{2 \varepsilon}{3} + 2 \alpha + \frac{2 \delta}{\eta} + \frac{4 \alpha^{\halb}}{\eta}\\
& \stackrel{(\ref{normh}),(\ref{barphitildephi})}{\le} & \|[q h \, , \, \tilde{\varphi} \iota \psi(b)] \|  + \frac{2 \varepsilon}{3} + 2 \alpha + \frac{2 \delta}{\eta} + \frac{4 \alpha^{\halb}}{\eta} + \frac{2 \delta}{\eta}\\
& \stackrel{(\ref{qhtildephi}),(\ref{constants})}{<} & \varepsilon \, .
\end{eqnarray*}
Similarly, we check that
\begin{eqnarray*}
\lefteqn{\|\be_{B} b \be_{B} - \theta \sigma \iota \psi(b)\|}\\
& \stackrel{(\ref{btheta}),(\ref{qsigma})}{\le} & \|\theta(q (b \otimes \be_{\Zh}) q - \sigma \iota \psi(b))\| + \frac{\varepsilon}{3} \\
& \stackrel{(\ref{cp-appr})}{\le} & \|q (\varphi \psi (b) \otimes \be_{\Zh}) q - \sigma \iota \psi(b)\| + \frac{\varepsilon}{3} + \alpha\\
& \stackrel{(\ref{normh}),(\ref{barphitildephi}),(\ref{qh})}{\le} & \|q h \bar{\varphi}(\be_{F}) (\varphi \psi (b) \otimes \be_{\Zh}) q - \sigma \iota \psi(b) \| + \frac{\varepsilon}{3} + \alpha + \frac{\delta}{\eta} \\
& \stackrel{(\ref{barphi})}{=} & \|q h (\sum_{i=0}^{n} \varphi(\be_{F_{i}}) \varphi \psi (b) \otimes \varrho(e_{i+1})) q - \sigma \iota \psi(b) \| + \frac{\varepsilon}{3} + \alpha + \frac{\delta}{\eta} \\
& \stackrel{\ref{multiplicative-domain},(\ref{cp-appr}),(\ref{barphi})}\le & \|q h \bar{\varphi} \psi (b)  q - \sigma \iota \psi(b) \| + \frac{\varepsilon}{3} + \alpha + \frac{\delta}{\eta} + \frac{2 \alpha^{\halb}}{\eta} \\
& \stackrel{(\ref{normh}),(\ref{barphitildephi})}{\le} & \|q h \tilde{\varphi} \iota \psi(b) q - \sigma \iota \psi(b) \| + \frac{\varepsilon}{3}+ \alpha + \frac{\delta}{\eta} + \frac{2 \alpha^{\halb}}{\eta} + \frac{\delta}{\eta} \\
& \stackrel{(\ref{sigma-factorization}),(\ref{qhtildephi}),(\ref{constants})}{<} &  \varepsilon \, .
\end{eqnarray*}
Finally, let $\tau \in T(A)$ be a tracial state on $A$. Then $\tau \circ \theta$ is a tracial state on $A \otimes \Zh$, so it is of the form 
\[
\tau \circ \theta = \tau' \otimes \bar{\tau}
\]
for some $\tau' \in T(A)$ (again, $\bar{\tau}$ denotes the unique tracial state on $\Zh$). We obtain
\begin{eqnarray*}
\tau(\be_{B}) & \stackrel{(\ref{Bdef})}{=} & \tau  \theta \sigma(\be_{\tilde{F}}) \\
& \stackrel{(\ref{sigmadef}),(\ref{qdef})}{\ge} & (\tau' \otimes \bar{\tau})(\tilde{\varphi}(\be_{\tilde{F}})) - \eta \\
& \stackrel{(\ref{barphitildephi})}{\ge} & (\tau' \otimes \bar{\tau})(\bar{\varphi}(\be_{F})) - \eta - \delta \\
& \stackrel{(\ref{barphi})}{=} & \sum_{i=0}^{n} \tau'(\varphi(\be_{F_{i}})) \cdot \bar{\tau}(\varrho(e_{i+1})) - \eta - \delta \\
& \stackrel{(\ref{bartaubeta})}{\ge} & \beta \cdot \sum_{i=0}^{n} \tau'(\varphi(\be_{F_{i}}))  - \eta - \delta \\
& = & \beta \cdot \tau'(\varphi(\be_{F})) - \eta - \delta \\
& \stackrel{(\ref{cp-appr})}{\ge} & \beta \cdot (1- \alpha)  - \eta - \delta \\
& \stackrel{(\ref{betamu}),(\ref{mu})}{>} & 2 \mu - \alpha - \eta - \delta \\
& \stackrel{(\ref{constants})}{>} & \mu \, .
\end{eqnarray*}
Therefore, properties (i), (ii) and (iii) of Lemma \ref{characterization-tr0} are satisfied.  Since $A$ has real rank zero and finite decomposition rank, it satisfies Blackadar's second fundamental comparability property by \cite{W4}, Theorem 3.9. Now by Lemma \ref{characterization-tr0}, $A$ has tracial rank zero.
\end{nproof}

\bn
The preceding theorem can be slightly rephrased as follows:

\begin{ncor}
If $A$ is a separable simple unital $C^{*}$-algebra with real rank zero and finite decomposition rank, then $A \otimes \Zh$ is tracially AF.
\end{ncor}

\begin{nproof}
This follows from Theorems \ref{rr0>wu} and \ref{Z-stable>TAF}.
\end{nproof}
\en

\section{Classification results}

In this last section we apply Lin's classification theorem for tracially AF algebras and results of Elliott (in the ASH case) and Dadarlat, Elliott and Gong (in the AH case) to derive some corollaries of Theorem \ref{Z-stable>TAF}.

\bn
In \cite{L2}, Lin has confirmed the Elliott conjecture for the class of simple unital tracially AF algebras which satisfy the UCT. As a consequence we have the

\begin{ncor}
Let $A$ and $B$ be separable simple unital $C^{*}$-algebras with real rank zero and finite decomposition rank; suppose $A$ and $B$ satisfy the UCT. Then, $A \otimes \Zh$ and $B \otimes \Zh$ are isomorphic iff their Elliott invariants are.
\end{ncor}
\en

\bn
\begin{ncor}
Let $A$ be a separable simple unital $C^{*}$-algebra with real rank zero and finite decomposition rank; suppose $A \otimes \Zh$ satisfies the UCT. Then:
\begin{itemize}
\item[(i)] $A \otimes \Zh$ is AH of topological dimension at most 3.
\item[(ii)] $A \otimes \Zh$ is ASH of topological dimension at most 2.
\item[(iii)] $\dr (A \otimes \Zh)$ is at most 2.
\item[(iv)] $A \otimes \Zh$ is approximately divisible. 
\item[(v)] $A$ is $\Zh$-stable iff $A$ is approximately divisible. 
\end{itemize}
\end{ncor}

\begin{nproof}
(i), (ii) and (iii) follow from results of Dadarlat, Elliott and Gong as in \cite{W4}, Corollary 6.4. By \cite{EGL}, an AH algebra of bounded topological dimension is approximately divisible. Conversely, an approximately divisible $C^{*}$-algebra is $\Zh$-stable by \cite{TW2}.
\end{nproof}
\en

\end{document}